\documentclass[11pt,leqno]{article}
\usepackage{amsmath,amssymb,latexsym}
\usepackage{amscd}
\usepackage{mathrsfs}
\usepackage{amsthm}
\usepackage{amsfonts}

\newcommand{\bproof}{{\raggedright\textbf{Proof.}} \ }
\newcommand{\proofmain}{{\raggedright\textbf{Proof of Theorem \ref{th_RIforvDlambda}.}} \ }
\newcommand{\bqed}{\hspace*{\fill} $\Box $\medskip}

\newcommand{\cK}{\mathcal{K}}
\newcommand{\cE}{\mathcal{E}}
\newcommand{\cD}{\mathcal{D}}
\newcommand{\cI}{\mathcal{I}}

\newcommand{\cB}{\mathcal{B}}

\newcommand{\C}{\mathbb{C}}
\newcommand{\N}{\mathbb{N}}
\newcommand{\R}{\mathbb{R}}

\newcommand{\cEK}{\mathcal E(K)}
\newcommand{\ern}{\mathcal E(\R^n)}
\newcommand{\cIK}{\mathcal{I}_K}
\newcommand{\rnrn}{:\R^n\to\R^n}
\newcommand{\ekek}{:\mathcal E(K)\to\mathcal E(K)}
\newcommand{\ernern}{:\mathcal E(\R^n)\to\mathcal E(\R^n)}

\newcommand{\clri}{continuous linear right inverse for }

\newcommand{\id}{\operatorname{id}\,}

\theoremstyle{plain}
\newtheorem{Th}{Theorem}[section]
\newtheorem{Cor}[Th]{Corollary}
\newtheorem{Lem}[Th]{Lemma}
\newtheorem{Prop}[Th]{Proposition}

\newtheorem*{Prob}{Problem}

\theoremstyle{definition}
\newtheorem{Def}[Th]{Definition}

\theoremstyle{remark}
\newtheorem{Exam}[Th]{Example}

\title{{\sc Right inverses for partial differential operators on spaces of Whitney functions}}
\author{{\sc Tomasz Cia{\'s}}}
\date{}

\begin{document}

\maketitle

\begin{abstract}
For $v\in\R^n$ let $K$ be a compact set in $\R^n$ containing a suitable smooth surface and such that the intersection $\{tv+x:t\in\R\}\cap K$ is 
a closed interval or a single point for all $x\in K$. 
We prove that every linear first order differential operator with constant coefficients in direction $v$ on space of Whitney functions $\cEK$ admits 
a continuous linear right inverse.  
\end{abstract}


\section{Introduction}\label{introduction}
In this paper we consider linear partial differential operators $P(D)$ with constant coefficients on the space of smooth Whitney functions 
$\cEK$ on a given compact set $K\subset\R^n$. By surjectivity of $P(D)$ on the space $\ern$ of smooth functions on $\R^n$ 
(see \cite[Cor. 3.5.2]{Hor}) it follows that $P(D)$ on $\cEK$ is surjective as well (see also \cite[p. 40]{Fr1}). 
In other words, for all $f\in\cEK$ the equation 
\[P(D)g=f\]
has a solution $g\in\cEK$. 
Now we can ask if it is possible to give solutions in the continuous and linear way.
More precisely, we are interested in the following problem: does $P(D)$ admit a continuous linear right inverse, i.e., an operator $S\ekek$ such 
that $P(D)\circ S=\id_{\cEK}$? So far, we know very little, and even there is no negative example.    

We say that a compact set $K\subset\R^n$ has the extension property if there exists a continuous linear extension operator $E:\cEK\to\ern$, i.e., 
$E$ satisfies the identity $r_K\circ E=\id_{\cEK}$, where $r_K:\cE(\R^n)\to\cE(K)$ is the continuous restriction operator (for the precise definition
see Section \ref{pre}). It is well known (see Prop. \ref{prop_K_ext_prop} below) that for $K$ with the extension property \emph{every} 
linear partial differential operator $P(D)$ with constant coefficients on $\cE(K)$ has a continuous linear right inverse. 
Compact sets with this property are very well characterized in terms of so called property $(DN)$ ($K$ has the extension
property if and only if $\cE(K)$ has the property $(DN)$; see \cite[Th. 3.3]{Fr2} or \cite[Folgerung 2.4]{Tid}) 
but the geometric characterization is still not known.

The case of compact sets without the extension property is much more complicated and to our best knowledge there is no (nontrivial) results in 
that case so far. In this paper we consider compact sets $K$ and partial differential operators $P(D)$ such that there exists a continuous linear right inverse 
$\widetilde S$ for $P(D)$ on $\ern$ so that $\widetilde S(\cI_K)\subseteq\cI_K$
($\cI_K$ stands for the ideal of functions flat on $K$).
Hence the operator 
\[Sf:=r_K(\widetilde{S} F)\] 
is defined independently of the choice of the extension $F\in\ern$ of $f\in\cE(K)$ and defines a continuous linear right inverse on $\cE(K)$ 
for the given differential operator (see Propostion \ref{prop_S0PD}). 
It appears that in the case of first order differential operator with constant coefficients in direction $v\in\R^n$ we can obtain 
such a $\widetilde{S}$ if a compact set $K$ contains a suitable smooth surface and it is so that the intersection 
$\{tv+x:t\in\R\}\cap K$ is a closed interval or a single point for all $x\in K$ ($v$-normal sets with a smooth surface defined in \ref{def_vnormalset}).
This is the main result of this paper (Theorem \ref{th_RIforvDlambda}). 

We divide the proof in a few steps. We start with the case of normal set in direction $e_j$ which contains the zero surface. 
Then using composition operators we pass to normal sets containing a smooth surface. 
The last step - right inverse in the case of $v$-normal sets with a smooth surface - is the result of a "rotation", i.e., 
the compostion with an appriopriate orthogonal linear map. 
Composing obtained in this way rigth inverses we get a right inverse in the case of compact sets which are normal in several directions simultaneously 
(Corollary \ref{cor_RI_v1vk}).

\section{Preliminaries}\label{pre}

Let us fix $n\in\N$ and let $\mathcal E(\R^n)$ denote the space of smooth functions on $\R^n$ with its natural Fr\'echet space topology. 
For an index $\alpha=(\alpha_1,\ldots,\alpha_n)\in\N_0^n$ we write $D^\alpha:=D_1^{\alpha_1}\ldots D_n^{\alpha_n}$, 
where $D_j\ernern$, $D_j:=\frac{\partial}{\partial x_j}$. 
More generally, for $P\in\C[x_1,\ldots,x_n]$ a polynomial of degree $N$, we consider 
the partial differential operator $P(D)\ernern$,
\[P(D):=\sum_{|\alpha|\leq N}\frac{D^{\alpha} P(0)}{\alpha!}D^\alpha,\]
where $|\alpha|:=\alpha_1+\ldots+\alpha_n$, $\alpha!:=\alpha_1!\cdot\ldots\cdot\alpha_n!$.
In particular, for a fixed $v=(v_1,\ldots,v_n)\in\R^n$, the directional derivative $\sum_{j=1}^nv_jD_j$ is denoted by $D_v$.

For a compact set $K\subset\R^n$ we define the restriction operator $r_K:\mathcal E(\R^n)\to\prod_{\beta\in\N_0^n}C(K)$,
\[r_KF:=((D^{\beta}F)\mid_{K})_{\beta\in\N_0^n}\] 
and by $\cI_{K}$ we denote the ideal of smooth functions which are flat on $K$, namely $\cI_K:=\ker r_K$.  
Let $\mathcal E(K)$ denote the space of Whitney functions on $K$,  
\[\mathcal E(K):=\{f=\big(f^{\beta}\big)_{\beta\in\N_0^n}:r_KF=f\quad\mathrm{for\phantom{x}some}\quad F\in\mathcal E(\R^n)\}.\]
The topology in $\mathcal E(K)$ is defined as the finest topology such that the restriction operator 
$r_K:\mathcal E(\R^n)\to\mathcal E(K)$ is continuous. It is easy to see that $\mathcal E(K)$ with this topology is a Fr\'echet space. 
For more information about spaces of Whitney functions, we refer to \cite{Fr2} and \cite{Mal}.

We introduce partial differential operators $D_j$, $D^\alpha$, $P(D)$ and $D_v$ on $\mathcal E(K)$ using the same notation as in the case of the space
$\mathcal E(\R^n)$. In this way we denote $D_jf:=(f^{\beta+e_j})_{\beta\in\N_0^n}$, $D^\alpha f:=(f^{\beta+\alpha})_{\beta\in\N_0^n}$, 
$P(D):=\sum_{|\alpha|\leq N}\frac{D^{\alpha} P(0)}{\alpha!}D^\alpha$ and $D_v:=\sum_{j=1}^nv_jD_j$, where $e_j$ is the vector in $\R^n$ which 
$j$-th coordinate equals 1 and the others equal 0. 

The following result is well known (for the proof see also \cite[Prop. 6.1]{Fr1}). 
\begin{Prop}\label{prop_K_ext_prop}
Let $K$ be a compact set in $\R^n$ with the extension property. Then every linear partial differential operator $P(D)$ with constant coefficients 
admits a continuous linear right inverse.
\end{Prop}
\bproof
Let $\cD'(\R^n)$ denote the space of distributions on $\R^n$ with compact support. 
If $\sigma\in\cD'(\R^n)$ is a fundamental solution for $P(D):\cD'(\R^n)\to\cD'(\R^n)$, $E:\cEK\to\cE(\R^n)$ is a linear continuous
extension operator and $\psi$ is a test function so that $\psi\equiv 1$ on a neighborhood of $K$, then 
\[Sf:=r_K(\sigma\ast(\psi\cdot Ef))\] 
is a \clri $P(D)$ on $\cEK$. 
\bqed

Let $||x||:=(\sum_{k=1}^nx_k^2)^{\frac{1}{2}}$ denote the euclidean norm of $x\in\R^n$. We denote by 
$\langle x,y\rangle:=\sum_{k=1}^nx_ky_k$ the scalar product of vectors $x,y\in\R^n$. For $v\in\R^n$ let
\[H_v^{(n)}:=\{x\in\R^n:\langle x,v\rangle=0\}\]
be the hyperplane in $\R^n$ orthogonal to $v$ containing 0. In particular, for $j=1,\ldots,n$ we write 
\[H_j^{(n)}:=\{x\in\R^n:\langle x,e_j\rangle=0\}.\] 
Let $K_v\subset H_v^{(n)}$ be compact and let $\phi_v,\psi_v:K_v\to\R$, $\phi_v\leq\psi_v$. We denote
\[\cK(K_v,\phi_v,\psi_v):=\{tv+x:x\in K_v, t\in[\phi_v(x),\psi_v(x)]\}.\] 

\begin{Def}\label{def_vnormalset}
Let $v_1,\ldots,v_k\in\R^n$.
We say that a set $K\subset\R^n$ is $(v_1,\ldots,v_k)$-normal if there exist compact sets $K_{v_m}\subset H^{(n)}_{v_m}$ and functions 
$\phi_{v_m},\psi_{v_m}:K_{v_m}\to\R$, $\phi_{v_m}\leq\psi_{v_m}$ such that $K=\cK(K_{v_m},\phi_{v_m},\psi_{v_m})$ for $m=1,\ldots,k$. 
Futhermore, if there exist $\gamma_{v_m}\in\mathcal E(K_{v_m})$ such that $\phi_{v_m}\leq \gamma_{v_m}^0\leq\psi_{v_m}$ for $m=1,\ldots,k$, then we 
say that $K$ is $(v_1,\ldots,v_k)$-normal with a smooth surface. In particular, if $\phi_{v_m}\leq 0\leq\psi_{v_m}$ for $m=1,\ldots,k$, 
then we say that $K$ is $(v_1,\ldots,v_k)$-normal with the zero surface. If $v_1=e_{l_1},\ldots,v_k=e_{l_k}$ for some $l_1,\ldots,l_k\in\{1,\ldots,n\}$, 
then we write for simplicity $(l_1,\ldots,l_k)$-normal instead of $(e_{l_1},\ldots,e_{l_k})$-normal. If $k=1$ we write $v$-normal ($j$-normal) 
for appropriate $v\in\R^n$ ($j\in\{1,\ldots,n\}$).  
\end{Def}

In view of Proposition \ref{prop_K_ext_prop}, the problem of the existence of a continuous linear right inverse for $P(D)$ on $\cE(K)$ is interesting
only for compact sets without the extension property. 
We give below examples of such sets which are simultaneously in the class described in Definition \ref{def_vnormalset}.

\begin{Exam}
(i) Let $K_1:=\{(x_1,x_2)\in\R^2:0\leq x_1\leq1,0\leq x_2\leq \exp(-1/x_1)\}$ and 
$K_2:=\{(x_1,x_2)\in\R^2:0\leq x_1\leq1, x_1^s\leq x_2\leq x_1^s+\exp(-1/x_1)\}$, where $s\geq1$ is  not rational. Then $K_1$, $K_2$ are compact, $(e_1,e_2)$-normal sets with 
a smooth surface ($K_1$ has even the zero surface) and 
they do not have the extension property (see e.g. \cite[Ex. 3.12 and Ex. 4.15]{Fr2}). \\
(ii) Let $f\colon\R^{n-1}\to\R$ be a smooth function and let $K:=\{(x_1,\ldots,x_{n-1},f(x_1,\ldots,x_{n-1}))\in\R^n: (x_1,\ldots,x_{n-1})\in K_0\}$,
where $K_0$ is an arbitrary compact set in $\R^{n-1}$.
Then one can easily check that the space $\cE(K)$ does not have a continuous norm 
(consider functions $(x_n-f(x_1,\ldots,x_{n-1}))^r$, $r\in\N$)
so, clearly, it does not have the property $(DN)$ (see \cite[p. 359]{MeV} for definition). 
Hence, from the Tidten's characterization (see \cite[Folgerung 2.4]{Tid}), $K$ does not have the extension property and, obviously, $K$ is a compact, 
$e_n$-normal set with a smooth surface. 

\end{Exam}

\section{Right inverse in the case of $v$-normal sets}

Let us formulate our main result.

\begin{Th}
\label{th_RIforvDlambda}
Let $\lambda\in\C$, $v=(v_1,\ldots,v_n)\in\R^n$, $v\neq0$ and let $K\subset\R^n$ be a compact, $v$-normal set with a smooth surface. 
Then the differential operator $\sum_{j=1}^nv_jD_j-\lambda\ekek$ admits a continuous linear right inverse.
\end{Th}

We can easily pass from the case of the first order differential operator to partial differential operators of any given order using 
the following lemma.

\begin{Lem}
\label{lem_PDPjD}
Let $P\in\C[x_1,\ldots,x_n]$, $P=P_1\cdot\ldots\cdot P_k$ for some polynomials $P_1,\ldots, P_k\in\C[x_1,\ldots,x_n]$ and let $K\in\R^n$ be compact. 
Then $P(D)\ekek$ has a continuous linear right inverse if and only if every $P_j(D)$ has a continuous linear right inverse.
\end{Lem}
\bproof
Let $S$ be a continuous linear right inverse for $P(D)$. Then 
\[P_j(D)\circ(P_1(D)\circ\ldots\circ P_{j-1}(D)\circ P_{j+1}(D)\circ\ldots\circ P_k(D)\circ S)=P(D)\circ S=\id_{\cEK},\]
hence $P_1(D)\circ\ldots\circ P_{j-1}(D)\circ P_{j+1}(D)\circ\ldots\circ P_k(D)\circ S$ is a continuous linear right inverse for $P_j(D)$. 

Conversly, if $S_j$ is a \clri $P_j(D)$, then $S_k\circ\ldots\circ S_1$ is a \clri $P(D)$.
\bqed

\begin{Cor}
\label{cor_RI_v1vk}
Let $P_1,\ldots,P_k$ be complex polynomials of one variable, $v_1,\ldots,v_k\in\R^n$ and let $K\subset\R^n$ be a compact, 
$(v_1,\ldots,v_k)$-normal with a smooth surface set. Then 
\[P(D)=P_1(D_{v_1})\circ\ldots\circ P_k(D_{v_k})\]
admits a continuous linear right inverse.
\end{Cor}
\bproof 
Follows from Theorem \ref{th_RIforvDlambda} and Lemma \ref{lem_PDPjD}.
\bqed

In the case of compact set whithout a smooth surface and without the extension property our method fails. 
It would be worth to solve the following problem.   
\begin{Prob}
Give an example of compact, 1-normal set $K$ in $\R^n$ which has no smooth surface, 
without the extension property and such that $D_1\ekek$ admits a continuous linear right inverse.  
\end{Prob}

In order to prove Theorem \ref{th_RIforvDlambda} we need several lemmas. First, we explain commutativity of differential operators with other types of 
operators. It is easy to prove the following lemma. 

\begin{Lem}
\label{lem_rKPD}
Let $K\subset\R^n$ be a compact set. Then $P(D)\circ r_K=r_K\circ P(D)$ for every polynomial $P\in\C[x_1,\ldots,x_n]$. 
\end{Lem}

Let us recall that for a smooth function $\Phi:\R^n\to\R^n$ the composition 
operator $\widetilde{C}_\Phi:\mathcal E(\R^n)\to\mathcal E(\R^n)$ is defined by the formula $\widetilde{C}_\Phi F=F\circ\Phi$.
\begin{Lem}
Let $\Phi:\R^n\to\R^n$ be a smooth function and let $K\subset\R^n$ a be compact set. 
Then $\widetilde{C}_\Phi\big(\mathcal{I}_{\Phi (K)}\big)\subset\mathcal I_K$. 
\end{Lem}
\bproof
Follows from the Fa\`a di Bruno formula (see e.g. \cite{Mis}).
\bqed\\
The preceding lemma allows us to define a composition operator $C_{\Phi,K}:\mathcal E(\Phi(K))\to\mathcal E(K)$, 
\[C_{\Phi,K}f:=r_K(\widetilde{C}_\Phi F),\] 
where $F$ is arbitrarily choosen function from $\mathcal E(\R^n)$ with $r_{\Phi(K)}F=f$. 
One can show that if $\Phi$ is a smooth bijection with the smooth inverse, then $C_{\Phi,K}$ is a continuous isomorphism with 
the inverse $C_{\Phi,K}^{-1}=C_{\Phi^{-1},\Phi(K)}$. 

\begin{Prop}
\label{prop_SPhiforPD}
Let $\Phi:\R^n\to\R^n$ be a smooth bijection with the smooth inverse and let $K\subset\R^n$ a be compact set. 
If $Q(D)\circ C_{\Phi^{-1},\Phi(K)}=C_{\Phi^{-1},\Phi(K)}\circ P(D)$ and $S:\mathcal E(K)\to\mathcal E(K)$ is a 
contionuous linear right inverse for $P(D):\mathcal E(\R^n)\to\mathcal E(\R^n)$, then $S_{\Phi}:\mathcal E(\Phi(K))\to\mathcal E(\Phi(K))$,
$S_{\Phi}:=C_{\Phi^{-1},\Phi(K)}\circ S\circ C_{\Phi,K}$
is a continuous linear right inverse for $Q(D):\mathcal E(\Phi(K))\to\mathcal E(\Phi(K))$.   
\end{Prop}
\bproof $S_{\Phi}$ is continuous as it is a composition of continuous operators. Clearly,
\begin{align*}
Q(D)\circ S_{\Phi}&= Q(D)\circ C_{\Phi^{-1},\Phi(K)}\circ S\circ C_{\Phi,K}=C_{\Phi^{-1},\Phi(K)}\circ P(D)\circ S\circ C_{\Phi,K}\\
&=C_{\Phi^{-1},\Phi(K)}\circ C_{\Phi,K}=\id_{\mathcal{E}(\Phi(K))}
\end{align*}
which means that $S_{\Phi}$ is a continuous linear right inverse for $Q(D)$. \bqed

\begin{Lem}
\label{lem_DjCPhi}
Let $j\in\{1,\ldots,n\}$, $\lambda\in\C$ be fixed, and let $\Psi=(\Psi_1,\ldots,\Psi_n):\R^n\to\R^n$ be a smooth bijection such that
\begin{displaymath}
D_j\Psi_l = \left\{ \begin{array}{ll}
1 & \textrm{if\quad $l=j$},\\
0 & \textrm{if\quad $l\neq j$}.\\
\end{array} \right.
\end{displaymath}
For a compact set $K\subset\R^n$ let $C_{\Psi}:=C_{\Psi,K}$ be a composition operator. Then
\[(D_j-\lambda)\circ C_{\Psi}=C_\Psi\circ(D_j-\lambda).\]
\end{Lem}
\bproof 
For $f\in\mathcal E(\Psi(K))$ let $F$ be its smooth extension, that is $r_{\Psi(K)}F=f$. By Lemma \ref{lem_rKPD}, we get
\begin{align*}
D_j(C_\Psi f)=& D_j(r_K(F\circ\Psi))
=r_K(D_j(F\circ\Psi))=r_K\bigg(\sum_{l=1}^n((D_lF)\circ\Psi)D_j\Psi_l\bigg)\\
=&r_K((D_jF)\circ\Psi)=C_{\Psi}(r_{\Psi(K)}(D_jF))
=C_{\Psi}(D_j(r_{\Psi(K)}F))=C_{\Psi}(D_jf),
\end{align*}
hence 
\[(D_j-\lambda)\circ C_{\Psi}=D_j\circ C_{\Psi}-\lambda C_{\Psi}=C_{\Psi}\circ D_j-C_\Psi\circ\lambda=C_\Psi\circ(D_j-\lambda).\]
\bqed

\begin{Lem}
\label{lem_DvCCDu}
Let $\lambda\in\C$, $u,v\in\R^n$, $u,v\neq 0$ and let $A\rnrn$ be a linear bijection such that $Au=v$. 
For a compact set $K\subset\R^n$ let $C_{A^{-1}}:=C_{A^{-1},K}$ be the composition operator. Then
\[(D_v-\lambda)\circ C_{A^{-1}}=C_{A^{-1}}\circ(D_u-\lambda).\] 
\end{Lem}
\bproof
Clearly, $\lambda\circ C_{A^{-1}}=C_{A^{-1}}\circ\lambda$. For $f\in\cE(A^{-1}(K))$ set $F\in\ern$ such that $r_{A^{-1}(K)}F=f$.
By Lemma \ref{lem_rKPD},
\[(D_v\circ C_{A^{-1}})(f)=D_v(C_{A^{-1}}f)=D_v(r_K(F\circ A^{-1}))=r_K(D_v(F\circ A^{-1}))\]
and, on the other hand,
\begin{align*}
(C_{A^{-1}}\circ D_u)(f)&=C_{A^{-1}}(D_uf)=C_{A^{-1}}(D_u(r_{A^{-1}(K)}F))=C_{A^{-1}}(r_{A^{-1}(K)}D_uF)\\
&=r_K(D_uF\circ A^{-1}), 
\end{align*}
hence it remains to observe that ($A^{-1}$ is linear) $D_v(F\circ A^{-1})=D_uF\circ A^{-1}.$
\bqed

Now we shall construct a right inverse on $\ern$. 
For $j\in\{1,\ldots,n\}$ and $\lambda\in\C$ we introduce the linear map $\widetilde{S}_{j,\lambda}:\mathcal E(\R^n)\to\mathcal E(\R^n)$, 
\[(\widetilde{S}_{j,\lambda}F)(x)=\int_0^{x_j}F(x^{(j,t)})e^{\lambda(x_j-t)}\mathrm{d}t,\]
where $x^{(j,t)}=(x_1,\ldots,x_{j-1},t,x_{j+1},\ldots,x_n)$.

\begin{Lem}
\label{lem_DalphatildeS}
Let $j\in\{1,\ldots,n\}$, $\lambda\in\C$ and $F\in\mathcal E(\R^n)$. Then
\[D^{\alpha}(\widetilde{S}_{j,\lambda}F)=\widetilde{S}_{j,\lambda}(D^\alpha F)\]  
for $\alpha$ with $\alpha_j=0$ and 
\[D_j^{\beta_j}D^{\alpha}(\widetilde{S}_{j,\lambda}F)
=\sum_{l=0}^{\beta_j-1}\lambda^lD_j^{\beta_j-l-1}D^{\alpha}F
+\lambda^{\beta_j}\widetilde{S}_{j,\lambda}(D^{\alpha}F)\]
for $\alpha$ with $\alpha_j=0$ and $\beta_j>0$.
\end{Lem}
\bproof
Let $l\neq j$. Then
\begin{align*}
(D_l(\widetilde{S}_{j,\lambda}F))(x)&=D_l\int_0^{x_j}F(x^{(j,t)})e^{\lambda(x_j-t)}\mathrm{d}t
=\int_0^{x_j}D_l\big(F(x^{(j,t)})e^{\lambda(x_j-t)}\big)\mathrm{d}t\\
&=\int_0^{x_j}(D_lF)(x^{(j,t)})e^{\lambda(x_j-t)}\mathrm{d}t,
\end{align*}
hence by induction we get the first formula. 

Now applying induction with respect to $\beta_j$ we will show that
\[D_j^{\beta_j}(\widetilde{S}_{j,\lambda}F)
=\sum_{l=0}^{\beta_j-1}\lambda^lD_j^{\beta_j-l-1}F+\lambda^{\beta_j}\widetilde{S}_{j,\lambda}F.\]
For $\beta_j=1$ we obtain
\begin{align*}
(D_j(\widetilde{S}_{j,\lambda}F))(x)&=D_j\bigg(\int_0^{x_j}F(x^{(j,t)})e^{-\lambda t}\mathrm{d}t\cdot e^{\lambda x_j}\bigg)\\
&=F(x)e^{-\lambda x_j}e^{\lambda x_j}+\lambda\int_0^{x_j}F(x^{(j,t)})e^{-\lambda t}\mathrm{d}t\cdot e^{\lambda x_j}
=F(x)+\lambda(\widetilde{S}_{j,\lambda}F)(x).
\end{align*}
Let us assume that
\[D_j^{\beta_j-1}(\widetilde{S}_{j,\lambda}F)
=\sum_{l=0}^{\beta_j-2}\lambda^lD_j^{\beta_j-l-2}F+\lambda^{\beta_j-1}\widetilde{S}_{j,\lambda}F.\]
Then
\begin{align*}
D_j^{\beta_j}(\widetilde{S}_{j,\lambda}F)=&D_j(D_j^{\beta_j-1}(\widetilde{S}_{j,\lambda}F))
=D_j\bigg(\sum_{l=0}^{\beta_j-2}\lambda^lD_j^{\beta_j-l-2}F+\lambda^{\beta_j-1}\widetilde{S}_{j,\lambda}F\bigg)\\
=&\sum_{l=0}^{\beta_j-2}\lambda^lD_j^{\beta_j-l-1}F+\lambda^{\beta_j-1}(F+\lambda\widetilde{S}_{j,\lambda}F)
=\sum_{l=0}^{\beta_j-1}\lambda^lD_j^{\beta_j-l-1}F+\lambda^{\beta_j}\widetilde{S}_{j,\lambda}F.
\end{align*}
Finally, we have
\[D_j^{\beta_j}D^{\alpha}(\widetilde{S}_{j,\lambda}F)=D_j^{\beta_j}(\widetilde{S}_{j,\lambda}(D^{\alpha}F))
=\sum_{l=0}^{\beta_j-1}\lambda^lD_j^{\beta_j-l-1}D^{\alpha}F
+\lambda^{\beta_j}\widetilde{S}_{j,\lambda}(D^{\alpha}F)\]
for $\alpha$ with $\alpha_j=0$ and $\beta_j>0$.
\bqed

\begin{Prop}
\label{prop_tildeSisRI}
Let $j\in\{1,\ldots,n\}$, $\lambda\in\C$ and let $K$ be a compact, $j$-normal set with the zero surface.
Then operator $\widetilde{S}_{j,\lambda}$ is a continuous linear right inverse for the differential operator 
$D_j-\lambda:\mathcal E(\R^n)\to\mathcal E(\R^n)$
such that $\widetilde{S}_{j,\lambda}(\mathcal I_K)\subset\mathcal I_K$.
\end{Prop}
\bproof 
Let $F\in\mathcal E(\R^n)$. By Lemma \ref{lem_DalphatildeS} (for $\beta_j=1$), we get
\[((D_j-\lambda)\circ \widetilde{S}_{j,\lambda})F=F\]
hence $\widetilde{S}_{j,\lambda}$ is a right inverse for $D_j-\lambda$. Now we shall show that $\widetilde{S}_{j,\lambda}$ is continuous. 
Let $(F_n)_{n\in\N}$ be a sequence in $\mathcal E(\R^n)$ such that 
\[\lim_{n\to\infty}F_n=F \quad \textrm{and} \quad \lim_{n\to\infty}\widetilde{S}_{j,\lambda} F_n=G\] 
for some $F,G\in\mathcal E(\R^n)$. This implies that $(F_n)_{n\in\N}$ is uniformly convergent on every compact subset of $\R^n$. 
Therefore, for fixed $x\in\R^n$ the sequence $(F_n(x^{(j,\cdot)})e^{\lambda(x_j-\cdot)})_{n\in\N}$ of smooth functions in one variable is uniformly 
convergent on the interval 
$[0,x_j]$ (or $[x_j,0]$ if $x_j<0$). Hence
\begin{align*}
(\widetilde{S}_{j,\lambda}F)(x)&=\int_0^{x_j}F(x^{(j,t)})e^{\lambda(x_j-t)}\mathrm{d}t
=\int_0^{x_j}\lim_{n\to\infty}F_n(x^{(j,t)})e^{\lambda(x_j-t)}\mathrm{d}t\\
&=\lim_{n\to\infty}\int_0^{x_j}F_n(x^{(j,t)})e^{\lambda(x_j-t)}\mathrm{d}t=\lim_{n\to\infty}(\widetilde{S}_{j,\lambda}F_n)(x)=G(x).
\end{align*} 
By the closed graph theorem, $\widetilde{S}_{j,\lambda}$ is continuous. Inclusion $\widetilde{S}_{j,\lambda}(\mathcal I_K)\subset\mathcal I_K$ follows 
immediately from the formulas given in Lemma \ref{lem_DalphatildeS}. 
\bqed

Now we get results which allow to transfer solutions from simpler cases to more complicated.
Firstly, let us remind a simple lemma about factorization of linear maps beetwen locally convex spaces by the quotient map 
(see e.g. \cite[Prop. 22.11]{MeV}).

\begin{Lem}
\label{lem_factorizationofT}
Let $T:X\to Z$ be a linear map beetwen locally convex spaces, $Y$ be a closed subspace of $X$ and $q:X\to X/Y$ be the quotient map. 
Let us also assume that $Y\subset \ker T$. Then there exists exactly one linear map $S:X/Y\to Z$ such that $T=S\circ q$. Moreover,
$S$ is continuous if and only if $T$ is continuous.
\end{Lem}

\begin{Prop}
\label{prop_S0PD}
Let $\widetilde{S}$ be a continuous linear right inverse for the differential operator $P(D):\mathcal E(\R^n)\to\mathcal E(\R^n)$
and let $K\subset\R^n$ be a compact set such that $\widetilde{S}(\mathcal I_K)\subset\mathcal I_K$. Let us define operator 
$S:\mathcal E(K)\to \mathcal E(K)$,
\[Sf:=r_K(\widetilde{S} F)\] 
where $F\in\mathcal E(\R^n)$, $r_K F=f$, is arbitrarily choosen. Then $S$ is a continous linear right inverse for the differential 
operator $P(D):\mathcal E(K)\to \mathcal E(K)$.
\end{Prop}
\bproof By Lemma \ref{lem_rKPD},
\[(P(D)\circ S)f=P(D)(r_K(\widetilde{S}F))=r_K(P(D)(\widetilde{S}F))=r_KF=f\]
hence $S$ is a right inverse for $P(D)$. From $\widetilde{S}(\mathcal I_K)\subset\mathcal I_K$ we get that $S$ is well defined and
its definition does not depend on the choice of the extension $F$ of $f$. Moreover,
$\cIK\subset\ker(r_K\circ\widetilde S)$ and, of course, we have $r_K\circ\widetilde{S}=S\circ r_K$. 
Now applying Lemma \ref{lem_factorizationofT} to the continuous operator $r_K\circ\widetilde{S}$ and the quotient map $r_K$ we obtain continuity of $S$.\bqed

\begin{Lem}
\label{lem_exK0Phi}
Let $j\in\{1,\ldots,n\}$ and let $K\subset\R^n$ be a compact, $j$-normal set with a smooth surface. 
Then there exist a compact, $j-normal$ set $K_0\subset\R^n$ with the zero surface and a smooth bijection $\Phi=(\Phi_1,\ldots,\Phi_n):\R^n\to\R^n$
with the smooth inverse such that 
\begin{displaymath}
D_j\Phi_l^{-1} = \left\{ \begin{array}{ll}
1 & \textrm{if\quad $l=j$},\\
0 & \textrm{if\quad $l\neq j$}\\
\end{array} \right.
\end{displaymath}
and $\Phi(K_0)=K$.
\end{Lem}
\bproof We have $K=\cK(K_j,\phi_j,\psi_j)$ for some real valued functions $\phi_j,\psi_j$ defined on a compact set $K_j$. 
Let $\gamma_j\in\mathcal E(K_j)$ satisfies $\phi_j\leq\gamma_j^0\leq\psi_j$ and set $\Gamma_j\in\ern$, $r_{K_j}\Gamma_j=\gamma_j$. 
Let us define $\Phi=(\Phi_1,\ldots,\Phi_n):\R^n\to\R^n$ by the formula
\[\Phi(x)=(x_1,\ldots,x_{j-1},x_j+\Gamma_j(x^{(j,0)}),x_{j+1},\ldots,x_n).\]
Clearly, $\Phi$ is smooth bijection with the smooth inverse and
\begin{displaymath}
\Phi_l^{-1}(x) = \left\{ \begin{array}{ll}
x_j-\Gamma_j(x^{(j,0)}) & \textrm{if\quad $l=j$},\\
x_l & \textrm{if\quad $l\neq j$},\\
\end{array} \right.
\end{displaymath}  
hence
\begin{displaymath}
D_j\Phi_l^{-1} = \left\{ \begin{array}{ll}
1 & \textrm{if\quad $l=j$},\\
0 & \textrm{if\quad $l\neq j$}.\\
\end{array} \right.
\end{displaymath}
Let $K_0:=\Phi^{-1}(K)$. It is easy to see that $K_0=\cK(K_j,\phi_j-\gamma_j^0,\psi_j-\gamma_j^0)$. 
Moreover, $\phi_j-\gamma_j^0\leq 0 \leq \psi_j-\gamma_j^0$, hence $K_0$ is a $j$-normal set with the zero surface. 
Finally, $K_0$ is compact as a continuous image of compact set $K$ and since $\Phi$ is a bijection we have $\Phi(K_0)=\Phi(\Phi^{-1}(K))=K$.
\bqed

\begin{Lem}
\label{lem_uvAK'K}
Let $u,v\in\R^n$, $||u||=||v||=1$ and let $K\subset\R^n$ be a compact, $v$-normal set with a smooth surface. 
Then there exist a compact, $u$-normal set $K'\subset\R^n$ with a smooth surface and a linear bijection $A\rnrn$ such that $Au=v$ and $A(K')=K$. 
\end{Lem}
\bproof 
By the Steinitz theorem and the Gram-Schmidt orthogonalization procedure, there is an orthogonal bijection $A\rnrn$, with $Au=v$ 
(i.e., $A^T=A^{-1}$, where $A^T$ is the conjugate operator for $A$). Then, it is easy to see that $A$ and $K':=A^{-1}(K)$
have desired properties. 
\bqed

Now we are ready to prove the main result.

\proofmain
Case of $v=e_j$ for some $j\in\{1,\ldots,n\}$, $K$ with the zero surface:
Combinig Propositions \ref{prop_tildeSisRI} and \ref{prop_S0PD} we get a right inverse $S_{j,\lambda}$ for $D_j-\lambda$ on $\mathcal E(K)$. 

Case of $v=e_j$ for some $j\in\{1,\ldots,n\}$, $K$ with an arbitrary smooth surface:
By Lemma \ref{lem_exK0Phi}, there exist a compact, $j$-normal set $K_0\subset\R^n$ with the zero surface and a smooth bijection 
$\Phi=(\Phi_1,\ldots,\Phi_n):\R^n\to\R^n$ with the smooth inverse such that 
\begin{displaymath}
D_j\Phi_l^{-1} = \left\{ \begin{array}{ll}
1 & \textrm{if\quad $l=j$},\\
0 & \textrm{if\quad $l\neq j$}\\
\end{array} \right.
\end{displaymath}
and $\Phi(K_0)=K$. From Lemma \ref{lem_DjCPhi} 
\[(D_j-\lambda)\circ C_{\Phi^{-1},K}=C_{\Phi^{-1},K}\circ(D_j-\lambda).\]
Now applying Proposition \ref{prop_SPhiforPD} to the function $\Phi$, the set $K_0$, the operator $D_j-\lambda$ and the operator 
$S_{j,\lambda}:\mathcal E(K_0)\to \mathcal E(K_0)$ from the previous case, we conclude that 
$C_{\Phi^{-1},K}\circ S_{j,\lambda}\circ C_{\Phi,K_0}$ is a \clri $D_j-\lambda\ekek$. 

Case of $||v||=1$, $K$ with an arbitrary smooth surface: By Lemma \ref{lem_uvAK'K}, there exist $e_1$-normal set $K'\subset\R^n$ 
with a smooth surface and a linear bijection $A\rnrn$ such that $Ae_1=v$ and $A(K')=K$. By Lemma \ref{lem_DvCCDu},
\[(D_v-\lambda)\circ C_{A^{-1}}=C_{A^{-1}}\circ(D_1-\lambda).\] 
From the previous case we obtain a continuous linear operator $S':\cE(K')\to\cE(K')$ such that $(D_1-\lambda)\circ S'=\id_{\cE(K')}$. 
Thus, by Proposition \ref{prop_SPhiforPD}, $S_A:=C_{A^{-1}}\circ S'\circ C_A$ is a continuous linear right inverse for 
$D_v-\lambda\ekek$. 

General case: Easily follows from the previous case.
\bqed

\flushleft{\textbf{Acknowledgements}. I wish to thank P. Doma\'nski and L. Frerick for several constructive remarks concerning this paper.

\vspace{1cm}
\begin{minipage}{7.5cm}
T. Cia\'s

Faculty of Mathematics and Comp. Sci.

A. Mickiewicz University in Pozna{\'n}

Umultowska 87

61-614 Pozna{\'n}, POLAND

e-mail: tcias@amu.edu.pl
\end{minipage}\

\end{document}